\documentclass[twocolumn,superscriptaddress,showpacs,amsmath,amssymb]{revtex4}

\usepackage{bm}

\newtheorem{thm}{Theorem}[section]
\newtheorem{prop}[thm]{Proposition}
\newtheorem{prop-def}[thm]{Proposition-Definition}
\newtheorem{lem}[thm]{Lemma}

\newtheorem{conj}[thm]{Conjecture}
\newtheorem{rem}[thm]{Remark}

\newtheorem{defn}[thm]{Definition}
\newtheorem{exmp}[thm]{Example}

\newcommand{\proofbegin}{\noindent{\it Proof:\,\,}}
\newcommand{\proofend}{\hfill$\Box$\bigskip}

\makeatletter
\def\@currentlabel{2.1}\label{e:dispaa}
\def\@currentlabel{2.21}\label{e:dispau}
\def\@currentlabel{2.22}\label{e:dispav}
\def\@currentlabel{2.23}\label{e:dispaw}
\def\@currentlabel{2.24}\label{e:dispax}
\def\theequation{\thesection.\@arabic\c@equation}
\makeatother

\makeatletter
\def\alphenumi{%
  \def\theenumi{\alph{enumi}}%
  \def\p@enumi{\theenumi}%
  \def\labelenumi{(\@alph\c@enumi)}}
\makeatother

\newcommand{\W}{\text{Dens}}

\newcommand{\Real}{\mathbb{R}}

\newcommand{\vect}{\mathfrak{X}}
\newcommand{\svect}{\mathfrak{X}_\mu}

\newcommand{\diff}{\mathcal{D}}
\newcommand{\sdiff}{\mathcal{D}_\mu}

\begin{document}

\title{Poisson geometry and first integrals \\of geostrophic equations}

\author{Boris Khesin}
\email{khesin@math.toronto.edu}
\affiliation{Department of Mathematics, University of Toronto, ON M5S 2E4, Canada}

\author{Paul Lee}
\email{plee@math.toronto.edu}
\affiliation{Department of Mathematics, University of Toronto, ON M5S 2E4, Canada}

\begin{abstract}
We describe first integrals of  geostrophic
equations, which are similar to the enstrophy invariants
of the Euler equation for an ideal incompressible fluid.
We explain the geometry behind this similarity,
give several equivalent definitions of the Poisson structure
on the space  of smooth densities on a symplectic manifold, and show
how it can be obtained via the Hamiltonian reduction from a
symplectic structure on the diffeomorphism group.
\end{abstract}



\maketitle

\section{Introduction}

The Euler equation of an ideal incompressible
$m$-dimensional fluid has a peculiar set of invariants: in
addition to the energy conservation, the Euler equation has the
helicity-type invariant in any odd $m$, and an infinite number 
of enstrophy-type invariants in any even $m$. Furthermore, 
these invariants are Casimir
functions, which implies that they are invariants of the Euler
equation for any choice of a Riemannian metric on the manifold
filled by the fluid.

In this paper we show how the same enstrophy-type
invariants appear in semi-geostrophic equations. These invariants are 
related to the Poisson geometry of the corresponding space of densities. 
Namely, for any Poisson manifold $M$ the space of densities on $M$
is also Poisson. The reason is that the space of functions on $M$ forms 
a Lie algebra  with respect to the Poisson bracket, while densities 
are objects dual to
functions, so their space forms a dual Lie algebra. Thus this dual space
gets equipped with the linear Kirillov-Kostant (or Lie-Poisson) structure, 
see \cite{Ambrosio, Lott, Wein}.
We show that this Poisson structure has several
equivalent descriptions and relate it to the symplectic geometry
of the diffeomorphism group of the manifold. In this paper we explore the role
of Casimirs and the corresponding group actions in
the Poisson geometry of these infinite-dimensional spaces. 

Recall that Poisson manifolds are foliated by symplectic leaves, and 
Casimir functions are functions constant on symplectic leaves. 
Equivalently, Casimir functions are those Hamiltonians which correspond 
to identically vanishing Hamiltonian vector fields on a Poisson manifold.
They are constants of motion for any Hamiltonian flow on the
manifold. For instance, the hydrodynamical Euler equation
on an odd-dimensional manifold has a helicity type Casimir, which generalises
the 3-dimensional helicity integral 
$$
I(v)=\int_M (v, {\rm curl}\,v)\,d\mu\,.
$$
For an even-dimensional manifold $M$ ($m=2n$) the Euler equation
has  an infinite number of enstrophy-type invariants:
$$
I_h(v)=\int_M h\left(\frac{({\rm curl}\,v)^n}{d\mu}\right)\,d\mu\,,
$$
where ${\rm curl}\,v$ is the vorticity 2-form for the velocity field $v$
on $2n$-dimensional manifold $M$ and $h$ is any function $\Real\to \Real$, 
see \cite{Serre, KhCh}. The latter integral
turns out to be similar to Casimir functions found on the space of 
densities on both even- and odd-dimensional manifold, as we discuss below.

One should mention that  the information on Casimir functions 
is useful for the study of the stability of Hamiltonian flows, e.g.
via the energy-Casimir method. Note that the existence of an infinite number 
of Casimirs for a given Hamiltonian system 
does not mean its complete integrability. These invariants merely single 
out the symplectic leaf  where the dynamics takes place, but do not specify
the dynamics along this leaf, cf. \cite{ArKh}.

 Below we also describe explicitly 
two natural Hamiltonian reductions leading to the Poisson 
structure on densities. A Hamiltonian reduction
is a two-step procedure for  reducing the dimension of a Hamiltonian system 
with symmetry: restriction to a given level set of first integrals, and 
taking the quotient along the symmetry group action. We prove that the Poisson
structure on densities can be obtained by the reduction 
from the symplectic structure considered in \cite{Do}.
More precisely,  let ${\mathcal Map}$ be the space of all maps 
from a manifold $M$ equipped with a volume form $\mu$
to a symplectic manifold $N$. The symplectic structure
on ${\mathcal Map}$ is given by averaging the pull-backs of the 
one on $N$ against the volume form $\mu$, see \cite{Do}. 
In Section \ref{ss-on-diff} we describe  what  
this general construction gives for diffeomorphism groups of symplectic 
manifolds and densities on them. 

\bigskip


\section{The Poisson structures on the density spaces and their Casimirs}

Let $M$ be a compact Poisson manifold with a Poisson bracket
$\{~,~\}$, and let $\W$ be the set of smooth volume forms on $M$
with total integral 1. (The  set  $\W$ can be given a smooth
topology and regarded as an infinite-dimensional smooth manifold.
See \cite{KrMi}. It is also a dense subset in the
$L^2$-Wasserstein space of Borel probability measures on $M$.) Any
smooth function $f$ on $M$ defines a linear functional on $\W$
whose value at a point $d\nu\in\W$ is given by the formula
$$
F_f(d\nu):=\int_M f\,d\nu\,.
$$
(For a non-compact $M$, e.g. for $M=\Real^n$, one can consider
functions $f$ with compact support.)

\begin{defn}\label{LottPoisson}
{\rm
Let $\mathcal P=\{F_f:\W\to\Real~|~f\in C^\infty(M)\}$
be the set of linear functionals $F_f$.
Define the  bracket on $\mathcal P$ by
$$
\{F_{f},F_{g}\}_\W(d\nu):= F_{\{f,g\}}(d\nu)
=\int_M \{f,g\}\,d\nu\,.
$$
}
\end{defn}

\begin{prop}\label{PoisWas}(see e.g. \cite{Wein, Lott})
\begin{enumerate}
   \item  The bracket $\{,\}_\W$ defines a Poisson structure on the
       density space $\W$.
  \item The symplectic
  leaves on $\W$ are orbits of the natural action of the group of Hamiltonian
  diffeomorphisms on  densities on $M$.
\end{enumerate}
\end{prop}

As we discussed in Introduction, the Poisson structure on the density space 
comes from the Poisson structure on the underlying manifold, as the Poisson-Lie
structure on the dual of the Lie algebra of Hamiltonian functions on $M$.
The statement (2) is proved in \cite{Lott}
for a symplectic $M$, but the proof extends verbatim
to a general Poisson manifold
$M$, provided that the group of Hamiltonian diffeomorphisms is understood as
that generated by flows of Hamiltonian fields on $M$.

\subsection{The symplectic case}

First consider  in more detail the case of a {\it symplectic} manifold $M$
of dimension $2n$. Let $\omega$ be a symplectic structure on $M$, which
generates the Poisson bracket $\{~,~\}$. Note that in this case the
Liouville form $\omega^n$ can be regarded as a natural choice
for a reference density $d\mu=\omega^n$.

\begin{prop}\label{prop:casimirs}
The Poisson bracket $\{,\}_\W$ admits infinitely many functionally
independent Casimirs. Namely, for any function
$h:\Real\to\Real$, the  functional
on $\W$ defined by
\[
C_h(d\nu):=\int_M h\Big(\frac{d\nu}{\omega^n}\Big)\,\omega^n
\]
is a Casimir, i.e. it is constant on symplectic leaves of this bracket in
the density space $\W$.
\end{prop}

\proofbegin
Since the  symplectic leaves of the Poisson structure $\{,\}_\W$
are orbits of the action of Hamiltonian flows on
the smooth Wasserstein space $\W$, it suffices to check that
the functions $C_h$ are invariant under this action. Now we have:

\begin{align*}
C_h(\phi^*d\nu)&
=\int_M h\Big(\frac{\phi^*d\nu}{\omega^n}\Big)\,\omega^n
\\&=\int_M h\Big(\frac{\phi^*d\nu}{\phi^*\omega^n}\Big)\,\phi^*\omega^n
\\&=C_h(d\nu)\,,
\end{align*}
where the last identity follows from the change of variable formula,
and the second one follows from conservation of $\omega$
under the Hamiltonian action: $\phi^*\omega=\omega$.
\proofend

Geometrically, these Casimirs capture all moments of the relative density
$d\nu$ with respect to the reference density $d\mu$.
The ratio function $\theta=d\nu/d\mu$ is preserved by
any Hamiltonian flow, and hence so are all its moments over the manifold $M$.

\begin{rem}
{\rm
Similar Casimirs arise in the case of the Euler equation
$$
\partial_t v + v\cdot\nabla v=-\nabla p
$$
for a divergence-free
vector field $v$ on any even-dimensional Riemannian manifold $M$
with volume form $d\mu$. Namely, one considers
the vorticity 2-form $du$ for the 1-form $u$ which is related to the vector
field $v$ by means of the metric on $M$. Then for any function
$h:\Real\to\Real$, the  functional on vorticities defined by
\[
I_h(du)=\int_M h\Big(\frac{(du)^n}{d\mu}\Big)\,d\mu
\]
is a Casimir for the action of diffeomorphisms  preserving the
``reference density''  $d\mu$, see \cite{Serre, KhCh} and Introduction. These
Casimirs also measure relative density of the generalized
vorticity $(du)^n$, which is frozen into the ideal fluid, with
respect to the volume form $d\mu$.}
\end{rem}

Conjecturally, a complete set of Casimirs is encoded in the
(Morse) graph with measure, associated to the function $\theta$ on
$M$. Its vertices correspond to critical points of $\theta$ on
$M$, and the edges correspond to pairs of critical points which
can be connected via nonsingular levels, while $d\theta$ defines
the measure on the graph. This construction has been used for
regular vorticity function in the 2D Euler equation (cf.
\cite{ArKh}), and is applicable to symplectic leaves in the
density space for any dimension.

\begin{exmp}
{\rm
Consider the following semi-geostrophic equation in (a domain of)
$\Real^2$:
$$
  \partial_t v_g + v\cdot\nabla v_g+Jv+\nabla f=0\,, \\
$$
where $J$ is the $90^\circ$-rotation operator on $\Real^2$,
$v$ is a divergence-free velocity field,
$v_g$ is the geostrophic velocity field ``defined by'' the relation
  $\nabla f=Jv_g$ for a potential $f$ in the domain, see \cite{Br}.
(This system is obtained from the two-dimensional
Euler equation in the rotating frame,
where we assume the Coriolis force to be constant in the domain,
and make the semigeostrophic approximation, see e.g. \cite{McOb}. )

Introduce the new potential $\tilde f(t,x):=|x|^2/2+f(t,x)$.
Consider the map $\phi_t(x)=\nabla \tilde f(t,\varphi_t(x))$,
where $\varphi_t$ is the flow of the divergence-free vector field
$v(t,\cdot)$ solving the above semi-geostrophic equation, and assume
that $\phi_t$ is a diffeomorphism for $t$ in some interval.
Then the family  $\phi_t$ descends to the following
Hamiltonian system on the density space by tracing how it pushes
the reference density $d\mu$. Namely, the form $d\nu_t:=
(\phi_t)_*d\mu$ satisfies the Hamiltonian system on the
space $\W$ with respect to the Poisson structure $\{,\}_\W$ and
the Hamiltonian $H^\W$ given by
\[
H^\W(d\nu)=-Wass^2(d\mu,d\nu)/2,
\]
where $Wass$ is the Wasserstein $L^2$-distance on $\W$.

The relative density $d\nu/d\mu$ discussed in Proposition
\ref{prop:casimirs} becomes
$$
\theta:=\frac{\phi_*d\mu}{d\mu}=\frac{(\nabla \tilde f)_*d\mu}{d\mu}
=\det(Hess~\tilde f)=\det(I+Hess~f)\,,
$$
where $Hess~ f $ is the Hessian matrix of the function $f$.
The latter expression for $\theta$ is known as the potential vorticity in
the semi-geostrophic equation, and is known to be frozen
into semi-geostrophic flow, similar to the standard vorticity of
an ideal two-dimensional fluid, see \cite{McOb}.

Thus Proposition \ref{prop:casimirs}
is a generalization of the Casimir property
of the potential vorticity to higher dimensions and to other Riemannian
metrics. Its frozenness property is shown to be related to the geometry of the
underlying Poisson structure $\{,\}_\W$
on the density space, rather than to the the specific Hamiltonian equation.
}
\end{exmp}

\subsection{The Poisson case}

Assume now that $M$ is a {\it Poisson} manifold whose symplectic leaves are
of codimension $\ge 1$, and $\lambda:M\to\Real$ is a smooth non-constant
Casimir function on $M$. It turns out that in this case
symplectic leaves of the Poisson bracket $\{,\}_\W$ still
have infinite codimension in $\W$, similar to the case of a symplectic $M$.

\begin{prop}
The  Poisson bracket $\{,\}_\W$ admits infinitely many functionally
independent Casimirs. Namely, for any function
$h:\Real\to\Real$, the  functional
\[
C_{h,\lambda}(d\nu):= \int_M (h\circ\lambda)\,d\nu
\]
is a Casimir on the density space $\W$.
\end{prop}

\proofbegin We check that the functionals $C_{h,\lambda}$ are
invariant under the Hamiltonian action:
\begin{align*}
C_{h,\lambda}(\phi^*d\nu)
&=\int_M  (h\circ\lambda)(x)\,\phi^*d\nu(x)
\\&=\int_M  (h\circ\lambda)(\phi(x))\,\phi^*d\nu(x)
\\&=C_{h,\lambda}(d\nu)\,,
\end{align*}
where we used the Casimir property of $\lambda$ on $M$:
$\lambda(\phi(x))=\lambda(x)$ for a Hamiltonian diffeomorphism $\phi$.
\proofend

Note that for symplectic leaves of codimension 1 on $M$,
one can think of invariants $C_{h,\lambda}$
as measuring the relative volume for the volume form $d\lambda\wedge\omega^n$
with respect to the reference density $d\mu$, where $\omega$ stands for
symplectic structure on the leaves in $M$.
This, in turn, is similar to the helicity-type  invariants for
the Euler equation on odd-dimensional manifolds, with the
important distinction, though, that for the
density space $\W$ one has not only one, but an infinite  number of Casimirs
regardless of the dimension of the manifold $M$.
\bigskip


\section{Symplectic structure on the diffeomorphism group of a
symplectic manifold}\label{ss-on-diff}

Let $(M,\omega)$ be a $2n$-dimensional symplectic manifold
and let $\diff$ be the space of all orientation
preserving diffeomorphisms of $M$. This is an infinite-dimensional
Lie group with the Lie algebra $\vect$ of all smooth vector fields on the
manifold $M$. The tangent space to the group $\diff$ at a
point $\phi$ consists of right translations of vector fields to $\phi$:
$T_\phi\diff=\{X\circ\phi~|~X\in\vect\}$. Fix the reference
volume form $d\mu=\omega^n$ on $M$.

\begin{defn}
{\rm
The diffeomorphism group $\diff$ can be equipped with the following
natural symplectic form $W^\diff$: given two tangent vectors
$X\circ\phi$ and $Y\circ\phi$ at $\phi\in\diff$ we set
\begin{align*}
W^\diff(X\circ\phi, Y\circ\phi)
:=\int_M \omega(X\circ\phi(x), Y\circ\phi(x))\,&d\mu(x)
\\=\int_M \omega(X, Y)\,(\phi^{-1})^*d\mu\,
=\int_M \omega(X, Y)\,\phi_*&d\mu\,.
\end{align*}
}
\end{defn}

As before, let $\W$ be the (smooth Wasserstein) space of all volume forms on
the manifold $M$ with total integral 1. The tangent space to this
infinite-dimensional manifold $\W$ at a point $d\nu$ consists of
smooth $2n$-forms on the manifold $M$ with zero integral. Denote the
tangent bundle of the smooth Wasserstein space by $T\W$.

Consider the natural projection
$\pi:\diff\to\W$ of diffeomorphisms into the volume forms on $M$,
according to how the diffeomorphisms move the reference density $d\mu$:
$\pi(\phi)=\phi_*\,(d\mu)$. This way the diffeomorphism group  $\diff$
can be regarded as the total space of the principal bundle over
the base  $\W$ with the structure group $\diff_\mu$ of all diffeomorphisms
preserving the volume form $d\mu$.

\begin{thm}\label{Poisson}
The symplectic structure $W^\diff$ on the
diffeomorphism group $\diff$ descends to the  Poisson structure
$\{~,~\}_\W$ on the  Wasserstein space $\W$.
\end{thm}

\proofbegin
The symplectic form $W^\diff$ is invariant under the $\sdiff$-action
of volume-preserving diffeomorphisms and hence under
the map $\pi$  it descends to a certain  Poisson structure
on the density space $\W$. We would like to show
that the corresponding quotient Poisson structure coincides with $\{~,~\}_\W$.

Let $f:M\to\Real$ be a function on the manifold $M$
and $F_f(d\nu)=\int_M f\,d\nu$ the corresponding linear functional on $\W$.
Consider the pullback $\bar F_f:=\pi^*F_f$ of this functional $F_f$
to the diffeomorphism group $\diff$ by the map $\pi$. Explicitly it is given by
$$
\bar F_f(\phi)=\int_M f\,d\nu
=\int_M f\,\phi_*(d\mu)
=\int_M (f\circ\phi)\,(d\mu)\,.
$$

Let $X_f$ be the Hamiltonian vector field for the Hamiltonian function $f$ on
the symplectic manifold $(M,\omega)$ and let $X_f^\diff$ be
the Hamiltonian vector field of the pullback functional $\bar F_f$ on
$(\diff,W^\diff)$, the diffeomorphism group $\diff$ equipped with the
symplectic structure $W^\diff$.

\begin{lem}\label{relation}
The Hamiltonian vector fields $X_f$ on $(M,\omega)$ and $X_f^\diff$ on
$(\diff,W^\diff)$ are related in the following way:
$$
X_f^\diff(\phi)=X_f\circ\phi.
$$
\end{lem}

\proofbegin
By the definition of the Hamiltonian field $X_{\bar F_f}$
at a point $\phi\in\diff$,
\[
W^\diff(X_f^\diff(\phi),Y\circ\phi)
=\langle d\bar F_f, Y\circ\phi\rangle
\]
for any vector field $Y\in\vect$.
On the other hand, by employing the definition of the pullback
and changing the variable, we rewrite
the latter expression as follows:
\[
\int_M \langle df_{\phi(x)}, Y\circ\phi(x)\rangle\,d\mu(x)
=\int_M \langle df_x, Y(x)\rangle\,\phi_*(d\mu)(x)\,.
\]
Now by the definition of the Hamiltonian field $X_f$ on $M$ this is equal to
\[
\int_M\omega(X_f(x),Y(x))\,\phi_*(d\mu)(x)
=\int_M\omega(X_f\circ\phi,Y\circ\phi)\,(d\mu)\,,
\]
which completes the proof of the lemma,
due to arbitrariness of the field $Y$.
\proofend

Returning to the proof of the theorem, we are going to compute
the Poisson bracket of the pullback functions $\bar F_{f}$ and
$\bar F_{g}$. By the definition, the value of the Poisson bracket
$\{,\}^\diff$, which is dual to the symplectic structure $W^\diff$ on the
diffeomorphism group, for these two functions is
\[
\{\bar F_{f},\bar F_{g}\}^\diff(\phi)
=W^\diff(X_{\bar F_{f}}(\phi),X_{\bar F_{g}}(\phi)).
\]
By using the lemma above and the change of variable,
the right-hand-side above becomes
\begin{align*}
\int_M \omega(X_{f}\circ\phi,X_{g}\circ\phi)\,d\mu
&=\int_M\omega(X_{f},X_{g}\,)\,\phi_*(d\mu)
\\&=\int_M \{f,g\}\,\phi_*(d\mu)
\\&=\int_M \{f,g\}\,d\nu\,,
\end{align*}
as required.
\proofend

\medskip
\begin{rem}
{\rm
This symplectic structure $W^\diff$ on the diffeomorphism group
can be viewed as a particular case of that considered  in \cite{Do}.
More generally, let $S$ be a compact manifold with a fixed volume form
$d\sigma$, while $(M,\omega)$ is a symplectic manifold.
The space ${\mathcal Map}$ of all maps $\rho:S\to M$ (of some
fixed homotopy class) has a natural symplectic structure. Namely,
the tangent space to ${\mathcal Map}$ at a point $\rho \in {\mathcal Map}$
is the space of sections of the bundle $\rho^*(TM)$ over $S$ and
the  symplectic structure is
$$
\Omega_f(v,w):=\int_M \rho^*\omega(v,w)\,d\sigma
$$
for a pair of sections $v,w$ of $\rho^*(TM)$. The group of
volume-preserving diffeomorphisms of $S$ defines a symplectic
group action on ${\mathcal Map}$. Donaldson  considers in \cite{Do}
the corresponding moment map and the Hamiltonian reduction of the
space ${\mathcal Map}$ under this group action. In our case, the two
manifolds $S$ and $M$ coincide, while the volume form $d\sigma$ is
the symplectic volume form $d\mu=\omega^n$. Then the
diffeomorphism group $\diff$ is an open  subset of $\mathcal
Map$ with the symplectic structure described above, and we consider
the action of the subgroup $\sdiff$ of volume-preserving
diffeomorphisms on it.
}
\end{rem}
\medskip
\begin{rem}
{\rm
The same Poisson structure on $\W$ was also defined in \cite{Ambrosio}
in slightly different terms, cf. \cite{Lott}.
For a symplectic manifold $(M,\omega)$ we fix a Riemannian metric
$\langle\,,\rangle$
and an almost complex structure $J$ compatible with the metric:
$\omega(u,v)=\langle u,Jv\rangle$.
Let $f$ be a function on the manifold $M$ and $\nabla f$ its gradient
with respect to the metric $\langle\, ,\rangle$.
The Hamiltonian field on $M$ for the Hamiltonian $f$ is $X_f=J\nabla f$.

Consider the distribution $\tau$ on the smooth Wasserstein space
defined at a point $d\nu\in\W$ by all possible
infinitesimal shifts of $d\nu$ by Hamiltonian fields:
$\tau_\nu:=\{L_{X_f}d\nu~|~ f\in C^\infty(M)\}$, where $L$ denotes
the Lie derivative along a vector field on $M$.
Define a  2-form on the distribution $\tau$ by
\[
\omega^\tau(L_{X_f}d\nu,L_{X_g}d\nu)=\int_M\omega(\nabla f,\nabla g)\,d\nu\,.
\]
In \cite{Ambrosio} it is shown that
the distribution $\tau_\nu$ is integrable
on the smooth Wasserstein space, and this 2-form is a well-defined
symplectic structure on the integral leaves of this distribution.
One can see that these leaves are exactly the symplectic leaves of the
Poisson structure $\{~,~\}_\W$ on the density space, while the symplectic
structure $\omega^\tau$ is dual to the Poisson structure discussed above:
\begin{align*}
\omega^\tau(L_{X_f}d\nu,L_{X_g}d\nu)
&=\int_M\omega(X_f,X_g)\,d\nu
\\&=\{F_f,F_g\}_\W(d\nu)\,.
\end{align*}
}
\end{rem}

\bigskip


\section{The two-dimensional case and geostrophic equations}

\subsection{The Noether theorem for an extra symmetry on the plane}

Return to the two-dimensional $M$ and
consider the smooth density space $\W$ for  $M=\Real^2$ with the
standard symplectic structure $\omega=dx_1\wedge dx_2$.
This induces the Poisson structure on $\W$, as described above.
There is the natural $SO(2)$-action by rotations on densities:
$d\nu\mapsto \varphi_*(d\nu)$,
where $d\nu\in\W$ is a measure and $\varphi\in SO(2)$.

Recall that for the standard measure $\omega$, the semi-geostrophic equation
is the Hamiltonian equation on $\W$ with the Hamiltonian function
$H(d\nu)=-Wass^2(\omega,d\nu)/2$, where $Wass$ is the Wasserstein distance on
densities.

\begin{prop}
The functional $K(d\nu):=\int_{\Real^2}|x|^2\,d\nu$ is a
first integral of the semi-geostrophic equation.
\end{prop}

\proofbegin
First we note that the $SO(2)$-action is Hamiltonian
with the Hamiltonian function given by the functional $K(d\nu)$ on $\W$.
Indeed, take the generator of the rotation group
with the Hamiltonian $\kappa(x)=|x|^2$ on $\Real^2$.
Then the corresponding action on densities in $\W$ is generated by
the field with Hamiltonian $K(d\nu):=\int_{\Real^2} \kappa\,d\nu
=\int_{\Real^2} |x|^2\,d\nu$,
while the corresponding action on diffeomorphisms in $\diff$ is
generated by the Hamiltonian $\bar K=\pi^*K$, cf. Lemma \ref{relation}.

Next, we see that the Wasserstein distance $H(d\nu)$ from any measure $d\nu$
to the standard measure $\omega$ is $SO(2)$-invariant, since so is $\omega$.
Thus the $SO(2)$-action is a
symmetry of the function $H(d\nu)=-Wass^2(\omega,d\nu)/2$, i.e.
the Hamiltonians  $H$ and $K$ are in involution on the density space $\W$
with respect to the Poisson structure $\{,\}_\W$.
In particular, $K(d\nu)$ is a conserved quantity for
the semi-geostrophic equation.
\proofend

This proposition naturally generalizes to any dimension:
If the Hamiltonian field with a Hamiltonian function $\kappa$ generates
an isometry of $M^{2n}$, and the reference density $d\mu=\omega^n$
is invariant  with respect to this isometry, then the Hamiltonian
field for $H(d\nu)=-Wass^2(d\mu,d\nu)/2$ on $\W$ has the first integral
$K(d\nu):=\int_{M} \kappa\,d\nu$.


\subsection{More Hamiltonian reductions to the density space}

Consider the case of a two-dimensional manifold $M$ in more detail.
In this section we would like to compare the symplectic geometry of the
diffeomorphism group $\diff(M)$ with that of the cotangent bundle
$T^*\sdiff(M)$ of the group of area-preserving diffeomorphisms of the surface
$M$.

As we discussed above, the group $\diff$ for an oriented surface (or,
for any symplectic manifold) $M$ can be equipped with a symplectic structure,
which descends to the Poisson structure on $\W$ under the projection
$\pi:\diff\to\W$, or, more precisely, under the Hamiltonian reduction
with respect to the $\sdiff$-action. Note that the space $\W$ is a convex
subset in the space $\Omega^2(M)$ of 2-forms on $M$:
$\W=\{d\nu\in\Omega^2(M)~|~d\nu>0,~\int_M d\nu=1\}$.

Now consider the group $\sdiff$ and its cotangent bundle $T^*\sdiff$.
Identify $T^*\sdiff\simeq \sdiff\times\svect^*$ by means of
right translations on the group. Note that the Lie algebra
$\svect$ consists of divergence-free vector fields on the surface $M$.
Such fields are described locally by their Hamiltonian (or, stream) functions.
First we assume that $M$ is a two-dimensional sphere $S^2$, so that the fields
are globally Hamiltonian. Then the algebra $\svect$
can be viewed as the Poisson algebra $C^\infty_0(M)$ of functions with
zero mean on $M$ (with respect to the reference density $d\mu=\omega$).
The corresponding (smooth)
dual space $\svect^*(M)=\Omega_0^2(M)$  consists of smooth  2-forms on $M$
with zero total integral.

By shifting this dual space to the reference density $d\mu$ one can regard
the (smooth) density space $\W$ as a convex subset in $\Omega_0^2(M)$:
$$
\W=\{d\mu+d\bar\nu~|~d\mu+d\bar\nu>0,~d\bar\nu\in \Omega_0^2(M)\}
\subset d\mu+\Omega_0^2(M)\,.
$$
After this shift to the reference density
$d\mu$ the diffeomorphism group $\diff\simeq \sdiff\times \W$ becomes
a subset in the cotangent bundle $T^*\sdiff\simeq \sdiff\times\Omega_0^2$.

Recall that the cotangent bundle  $T^*\sdiff$ has a natural symplectic
structure (denoted later by $W^{T^*}$),
which descends to the Poisson-Lie structure
on the dual Lie algebra $\svect^*\simeq \Omega_0^2$.
The latter is exactly the Poisson  structure $\{~,~\}_\W$
upon restriction to the density space $\W\subset\Omega_0^2$.

\begin{conj}
The natural symplectic structure $W^{T^*}$ on the cotangent bundle
 $T^*\sdiff$ coincides with the symplectic structure $W^\diff$
on the diffeomorphism group $\diff$, understood as a subset of
 $T^*\sdiff$ via the identification described above.
\end{conj}

In other words, not only coincide the Poisson structures on the density space
$\W$ understood by itself or as a part of the dual $\svect^*$,
but   presumably so do the corresponding symplectic structures before the
Hamiltonian reduction.
We note that the convex subset $\W$ of positive densities
is preserved under the diffeomorphism action on $\Omega_0^2$,
so the Poisson structure on $\Omega_0^2$ can indeed
be restricted to this subset.
\smallskip

In the case of a general surface $M$, divergence-free fields on $M$
may have multivalued Hamiltonians: $\svect(M)=C^\infty_0(M)\oplus H_1(M)$.
Respectively, the dual space $\svect^*(M)$ is a finite-dimensional
extension of $\Omega_0^2(M)$, since $\sdiff^*(M)\simeq
\Omega^1(M)/d\Omega^0(M)\simeq\Omega_0^2(M)\oplus H^1(M)$.
Now the density space $\W(M)$ can  be understood as a convex subset
in a plane of finite codimension in the dual space $\svect^*(M)$.

\smallskip

Note also that for a higher-dimensional symplectic $M$,
there is a natural map from the dual
$\svect^*\simeq \Omega^1(M)/d\Omega^0(M)\simeq\Omega_0^2(M)\oplus H^1(M)$
to the  space $\Omega^{2n}(M)$ of $2n$-forms:
$\rho: [u]\mapsto (du)^n$, where $[u]\in \svect^*$
is a 1-form $u$ modulo addition of an exact 1-form on $M$.
This map commutes with the natural $\sdiff$-action of
volume-preserving diffeomorphisms on forms, which explains the common origin
of the Casimirs on the space $\W$ of volume forms and
on the dual space $\svect^*$.

\medskip

Finally, we note that the Euler equation of an ideal fluid on the
two-dimensional $M$ is the Hamiltonian equation on $\Omega_0^2$
with respect to the Poisson-Lie structure
 $\{~,~\}_\W$ whose Hamiltonian function is the {\it energy} quadratic form on
 $\Omega_0^2$. It is interesting to compare this with (the projection of)
the semi-geostrophic equation, where the Hamiltonian function is
$H^\W(d\nu)=-Wass^2(d\mu,d\nu)/2$, the square
of the Wasserstein distance on $\W\subset d\mu+\Omega_0^2$.
This shift of the quadratic form to the reference density $d\mu$
is similar to the shift observed in
the infinite conductivity equation, and in the $f$-plane and $\beta$-plane
geostrophic equations, see \cite{ArKh, Holm-Zeitlin, Zeitlin}.

Before the Hamiltonian reduction, the semi-geostrophic equation
is a Hamiltonian system on the diffeomorphism group
with respect to the symplectic structure
$W^\diff$ and the Hamiltonian function $H^\diff$ evaluating how far
our diffeomorphism is from being  area-preserving:
\[
H^\diff(\phi)=\frac{1}{2}dist^2(\phi,\diff_\mu),
\]
where $dist$ is the distance in the ``flat'' $L^2$-type metric
on the  diffeomorphism group $\diff$ from $\phi\in \diff$ to
the subgroup $\diff_\mu$ of diffeomorphisms preserving the standard area form
$\omega=d\mu$ on $\Real^2$, see \cite{Br}.
This Hamiltonian is invariant under the action of the group $\sdiff$ of
area-preserving diffeomorphisms, and so it descends to the above
semi-geostrophic Hamiltonian system on the density space $\W$.

\bigskip





\begin{acknowledgments}
We are grateful to W.\,Gangbo, D.\,Goldman, and R.\,McCann
for useful discussions.
We are also indebted to the organizers of the conference ``Euler Equations: 250
years on'' in Aussois and  the workshop
``Optimal Transportation, and Applications to Geophysics and Geometry''
in Edinburgh for stimulating working atmosphere and hospitality.
This research was partially supported by an NSERC research grant.
\end{acknowledgments}


\end{document}